\documentclass[12pt,a4paper,dvipsnames]{amsart}

\usepackage[english]{babel}
\usepackage[T1]{fontenc}

\usepackage{etoolbox}
\apptocmd{\sloppy}{\hbadness 10000\relax}{}{}
\apptocmd{\sloppy}{\vbadness 10000\relax}{}{}

\usepackage[a4paper,top=2cm,bottom=2cm,left=1.5cm,right=1.5cm,marginparwidth=1.75cm]{geometry}

\usepackage{amsmath}
\usepackage{graphicx}
\usepackage[colorinlistoftodos]{todonotes}
\usepackage[pagebackref,hypertexnames=false, colorlinks, citecolor=red,linkcolor=blue, urlcolor=red]{hyperref}
\usepackage{diagbox,multirow} 
\usepackage{relsize}
\usepackage{comment}  

\usepackage{lmodern}
\usepackage[english]{babel}
\usepackage{amsmath}
\usepackage{amsfonts}
\usepackage{colonequals} 
\usepackage{amssymb}
\usepackage{amsthm}
\usepackage{stmaryrd}
\usepackage[all]{xy}

\usepackage{tikz}
\usepackage{float}
\usepackage{multicol}

\newtheorem{thm}{Theorem}
\newtheorem{rmq}{Remark}
\newcommand{\imaginaire}{\mathtt{i}}
\newcommand{\jj}{\mathtt{j}}
\newcommand{\kk}{\mathtt{k}}
\newcommand{\qU}{\mathfrak U}

\newcommand{\HH}{\mathbb H}

\begin{document}

\title[Limits of Stochastic Processes]{On some high-dimensional limits of matricial stochastic processes seen from a quantum probability perspective}
\author{Michaël Ulrich }
\thanks{Laboratoire de Mathématiques de Besançon\\ 16, route de Gray\\ 25000 Besançon, France\\ michael.ulrich@univ-fcomte.fr}

\date{February 2022}
\maketitle
\begin{abstract}
We generalize the result of block-wise convergence of the Brownian motion on the unitary group $U(nm)$ towards a quantum Lévy process on the unitary dual group $U\langle n\rangle$ when $m\rightarrow\infty$ obtained in \cite{Ulrich2015} by showing that the Brownian motions on the orthogonal group $O(nm)$ and the symplectic group $Sp(nm)$ also converge block-wise to this same quantum Lévy process.
\end{abstract}
\tableofcontents
\section*{Introduction}
Noncommutative mathematics' idea is to take a certain class of mathematical spaces $E$ whose properties can be adequately described by a space of functions over $E$, for instance $C(E)$ or $L^\infty(E)$, etc. Then, because this space of funtion has the structure of an algebra, possibly with some more properties, we replace the study of the class of spaces by the study of algebras having the same properties as the space of functions over $E$ but without assuming commutativity. Hence, noncommutative (or quantum) probability replaces the study of probabilized spaces $(\Omega,\mathcal F,\mathbb P)$ with the stdy of $*$-algebras $A$ endowed with a positive linear functional $\phi$ on $A$ such that $\phi(1)=1$. In the classical case, of course, $\phi$ is taken to be the expectation.\\
To generalize the notion of groups in this noncommutative setting, two approaches are possible.The one usually used is the quantum group approach, see \cite{Timmermann,KS,majid_1995,Kassel,CP}. Especially, a compact quantum group is defined to be a certain $C^*$-algebra denoted $C(G)$ equipped with a comultiplication  $\Delta:C(G)\rightarrow C(G)\otimes C(G)$ such that $(\Delta\otimes Id)\circ\Delta=(Id\otimes\Delta)\circ\Delta$ (coassociativity) and with another property called the quantum cancellation property. The second approach possible, the dual groups one, first initiated by Voiculescu in \cite{Voiculescu87}, consists in taking a $C^*$-algebra $A$ (actually we can take a mere $*$-algebra, and this is what we will do here) endowed with a comultiplication$\Delta:A\rightarrow A\sqcup A$ where $\sqcup$ is the free product in the category of $*$-algebras. This comultiplication needs to fulfill a similar coassociativity axiom.\\
Dual groups have been studied somewhat, though less than quantum groups. For instance, one can refer to \cite{Voiculescu87,FranzHabilitation,UlrichPhD,CU2016,kummerer1992,Ulrich2015,baraq2019,baraq2019b}. In particular, in \cite{Ulrich2015}, the author extended a previous result due to Biane in \cite{Biane1997}. Biane had shown that the Brownian motion on the unitary group $U(m)$ tends towards a multiplicative unitary Brownian motion when $m$ goes to infinity. Extending this, \cite{Ulrich2015} showed the Brownian motion on $U(mn)$ when $n$ is fixed and when the Brownian motion is seen block-wise, tends towards a Lévy process on the dual unitary group $U\langle n\rangle$. This Lévy process was nice in the sense that it was shown to be Gaussian, and, moreover, \cite{CU2016} showed that this quantum Lévy process tends towards the so-called Haar trace when $t$ goes towards infinity. So, it was a good candidate to be called a Brownian motion on $U\langle n\rangle$.\\
The question remained about what happens for the Brownian motion on the orthogonal group $O(nm)$ and on the symplectic group $Sp(nm)$ when they are seen block-wise and when $m$ goes towards infinity. We will show that they do converge towards the same quantum Lévy process that was found in \cite{Ulrich2015}. This is consistent with \cite[Theorem 2.2]{Levy2011}. While doing this, we will also introduce briefly the dual orthogonal group $O\langle n\rangle$ and the dual symplectic group $Sp\langle n\rangle$.

\section*{In honorem}
I dedicate this paper to Michael Schürmann on the occasion of his retirement. He was one of my PhD advisors. One our first encounter in 2014 in Greifswald, I had the occasion to show him the results I had obtained and that were to be published in \cite{Ulrich2015}. He was of course really supportive and happy about these results. We discussed together about the future and it was clear that I needed to investigate also the orthogonal and the symplectic cases. But, as often in research, other questions showed up and, trying to solve these latter, I put aside the orthogonal and symplectic question for a day when I would have more time... As is well know, such a day never exists in Academia! At the occasion of his retirement, I believe that it is fitting for me to finally complete this work, pushed aside for too long, as a token of my appreciation and friendship, as well as thankfulness for his mentorship.

\section{Dual groups and quantum probability}
It is not our goal here to explain in details the theory of dual (semi)groups and/or quantum probability. We refer to \cite{Ulrich2015, UlrichPhD, CU2016} for more information on this. Thus, let us just recall that a dual semigroup is a triple $(A,\Delta,\delta)$ with $A$ a $*$-algebra, and $\Delta:A\rightarrow\Delta\sqcup A$ and $\delta:A\rightarrow\mathbb C$ two $*$-homomorphisms, where $\sqcup$ denotes the free product of $*$-algebras. These $*$-homomorphisms must further fulfill the following relations:
$$\left(Id_A\bigsqcup \Delta\right)\circ\Delta=\left(\Delta\bigsqcup Id_A\right)\circ\Delta$$
and
$$\left(Id_A\bigsqcup \delta\right)\circ\Delta=Id_A=\left(\delta\bigsqcup Id_A\right)\circ\Delta$$
The dual semigroups we will consider here are:\\
1. The free dual semigroup $\mathcal M\langle n\rangle$, that is the dual semigroup defined on the free algebra $M_n$, defined on the complex field, and generated by $n^2$ generators $u_{i,j},1\leq i,j\leq n$, with no relations between them, and with the following $*$-homomorphims:
\begin{multicols}{2}
$$\Delta(u_{ij})=\sum_{k=1}^nu_{ik}^{(1)}u_{kj}^{(2)}$$
$$\delta(u_{ij})=\delta_{ij}$$
\end{multicols}
where $\delta_{ij}$ is Kronecker's symbol.\\
2. The unitary dual semi group $U\langle n\rangle$ defined on the $*$-algebra $U_n$ on the complex field, generated by the generators $u_{ij},1\leq i\leq j\leq n$ with the only relations given by $\sum_{k=1}^n u_{ki}^*u_{kj}=\delta_{ij}$ and $\sum_{k=1}^nu_{ik}u_{jk}^*=\delta_{ij}$, and with the $*$-homomorphisms:
\begin{multicols}{2}
$$\Delta(u_{ij})=\sum_{k=1}^nu_{ik}^{(1)}u_{kj}^{(2)}$$
$$\delta(u_{ij})=\delta_{ij}$$
\end{multicols}
It needs to be observed that, were we to add the additional commutation relation between the generators, we would then obtain the algebra generated by the maps $\tilde u_{ij}:U(n)\rightarrow \mathbb C$ defined on the (classical) unitary group by $\tilde u_{ij}(U)=U_{ij}$, ie associating the $(i,j)$ coefficient to a unitary matrix. In particular, the relations defining the $*$-algebra exactly express the fact that the matrix $\hat U$ containing the generators $u_{ij}$ as coefficients verifies $U^*U=I_n=UU^*$. This explains why $U\langle n\rangle$ is called the unitary dual semigroup. \\
3. We want to have an equivalent to the (classical) orthogonal group. However, we cannot just take the algebra generated by the $u_{ij}$ with the relations  $\sum_{k=1}^n u_{ki}u_{kj}=\delta_{ij}$ and $\sum_{k=1}^nu_{ik}u_{jk}=\delta_{ij}$, which would correspond to the intuitive orthogonal matricial relation $\,^tOO=I_n=O^tO$. Indeed, these relations do not behave well with regard to the coproduct $\Delta$, contrary to the unitary relations. So, we follow the lead of Voiculescu, who in \cite{Voiculescu87} decided to define the orthogonal dual semigroup $O\langle n\rangle$ as defined by the same algebra $U_n$ as for $U\langle n\rangle$, but this time taken over the field of real numbers. This corresponds to the fact that an orthogonal matrix verifies also the unitary relations (as all coefficients are real and thus not modified under conjugation) but that all coefficients being real, it is most natural to see the algebra of coefficient-maps as being over $\mathbb R$, rather that over $\mathbb C$.\\
4. The symplectic dual semigroup $Sp\langle n\rangle$. We will define it in the section dedicated to it, as we need to recall more details about the classical symplectic group $Sp(n)$.\\
Let us assume that we have two Lie groups $G$ and $G'$ such that $G\subset G'$. Then, if we have a map $f:G'\rightarrow \mathbb C$, then we also have a map $\tilde f:G\rightarrow \mathbb C$ where $\tilde f=f\circ j$ with $j$ the canonical embedding from $G$ to $G'$. Thus, we have a $*$-homomorphism $\tilde j:C(G')\rightarrow C(G)$, where $C(G)$ is the algebra of continuous functions defined on $G$, and the same for $C(G')$. As we can morally see the algebra of a dual semigroup as an algebra of noncommutative functions defined over the dual semigroup, we can decide to say that given two dual semigroups $\mathcal G$ and $\mathcal G'$, we say that $\mathcal G$ is contained in $\mathcal G'$ if there is a $*$-homomorphism $j:C(\mathcal G')\rightarrow C(\mathcal G)$ where $C(\mathcal G)$ designates the algebra of the dual semigroup $\mathcal G$. With this vocabulary, it is easy to see that $U\langle n\rangle$ is a dual subsemigroup of $\mathcal M\langle n\rangle$. Indeed, it suffices to take the trivial map, mapping any $u_{ij}$ of $M_n$ to the $u_{ij}$ of $U_n$.
\begin{rmq}
Let us remark that all of the aforementioned dual semigroups are actually dual groups, that is there exists a map $\Sigma:\mathcal G\rightarrow \mathcal G$ (where $\mathcal G$ is the algebra underlying the dual group) such that $(Id\sqcup \Sigma)\circ\Delta=1\circ\delta=(\Sigma\sqcup Id)\circ\Delta$. Indeed, in the three examples mentioned we can take this map, called the coïnverse, to be $\Sigma:u_{ij}\mapsto u_{ji}^*$. Nevertheless, to talk about Lévy processes, the structure of dual semigroup is enough. 
\end{rmq}
We also recall that a noncommutative unital $*$-probability space is given by a $*$-algebra $A$ and a tracial positive linear form $\phi$ on $A$ such that $\phi(1_A)=1$. Given a dual semigroup $\mathcal G$ whose algebra is denoted by $C(\mathcal G)$, we define a (noncommutative) random variable on $\mathcal G$ as the giving of a $*$-homomorphism $X:C(\mathcal G)\rightarrow A$. We observe that the direction of the map is the reverse of the usual direction in classical probability. This is consistent with the fact that we reason with dual structures. In the setting of quantum probability, there are five different notions of independence of random variables that can be used: classical (also called tensor) independence, freeness, boolean, monotone and antimonotone independence.\\
In all the sequel of the article, one $*$-probability space will be fixed.\\
Finally, let us recall the definition of a (quantum) Lévy process. Given a notion of independence $T$, it is a family $(j_{st})$ of quantum random variables such that
\begin{enumerate}
\item $(j_{rs}\sqcup j_{st})\circ\Delta=j_{rt}$ for all $0\leq r\leq s\leq t$ (increment property).
\item The quantum variables $j_{r_1s_1},\ldots,j_{r_ns_n}$ are independent in the sense of $T$, whenever $0\leq r_1\leq s_2<r_2\leq s_2\leq\ldots$
\item The distribution of $j_{st}$, ie the quantity $\phi\circ j_{st}$ depends only on $t-s$ (stationarity).
\item The quantum random variable $j_{st}$ converges towards $j_{ss}$ in distribution when $t\rightarrow s$ (weak continuity).
\end{enumerate}
More details can be obtained in \cite{FranzHabilitation}.Let us just observe that if $(X_t)_{t\geq0}$ is a (classical) Lévy process on a Lie group $G$, then, the process 
$$\begin{array}{ccc}j_{st}&:&C(G)\rightarrow L^{\infty-}(\Omega)\\&&f\mapsto f\circ X_{t-s}\end{array}$$
is a (quantum) tensor-independent Lévy process.
\begin{rmq}
If we have a dual group (instead of a mere dual semigroup), the only advantage, is that we can restrict ourselves to the study of $j_t=j_{0t}$.
\end{rmq}

\section{The quantum Lévy process}
\subsection{Definition of the process}
In the sequel, we will be interested in the process $\mathfrak U(t)$ defined as follows. We denote by $n$ a natural integer. We have a matrix $\mathfrak X(t)$ such that:
\begin{itemize}
\item[$\bullet$] At fixed $t$, the family $(\mathfrak X_{ij}(t))_{1\leq i\leq j \leq n}$ is free.
\item[$\bullet$] For any $i,j$, we have $\mathfrak X_{ij}^*=\mathfrak X_{ji}$.
\item[$\bullet$] For any $i$, $\mathfrak X_{ii}$ is a free additive Brownian motion
\item[$\bullet$] For any $i\neq j$, $\mathfrak X_{ij}=\mathfrak X_{ij}^{(1)}+\mathtt{i}\mathfrak X_{ij}^{(2)}$, where $\imaginaire$ is a square root of $-1$ and $\sqrt 2\mathfrak X_{ij}^{(p)}$, for $p\in\{1,2\}$, are two free additive Brownian motion, free with one another.
\end{itemize}
We remark that $\mathfrak X(t)$ is the free analogue of a Brownian motion on hermitian matrices and was already defined in \cite[Theorem 1][{Ulrich2015}. We now define the process $\mathfrak U(t)$ in a way that is similar to \cite[Theorem 1]{Ulrich2015}, but with a different renormalization, and by taking the left variant:
$$d\mathfrak U_t=\imaginaire \mathfrak U_t d\mathfrak X_t+\frac n2\mathfrak U_tdt$$
We follow \cite{kummerer1992} for the method of stochastic integration. We notice that we have $d\mathfrak X_{ij}.d\mathfrak X_{kl}=\delta_{il}\delta_{jk}dt$. Moreover, we have
$$d\mathfrak U_{ij}=\imaginaire \sum_{k=1}^n\mathfrak U_{ik}d\mathfrak X_{kj}+\frac n2\mathfrak U_{ij}dt$$
and
$$d\mathfrak U^*_{ij}=-\imaginaire \sum_{k=1}^nd\mathfrak X_{jk}\mathfrak U^*_{ik}+\frac n2\mathfrak U^*_{ij}dt$$
We then remember that the stochastic Itô integration by part formula states that:
\begin{eqnarray*}d(A_{i_1j_1}\ldots A_{i_nj_n})&=&\sum_{k=1}^n A_{i_1j_1}\ldots d(A_{i_kj_k})\ldots A_{i_nj_n}\\&+&\sum_{1\leq k<l\leq n}A_{i_1j_1}\ldots A_{i_{k-1}j_{k-1}}\phi(A_{i_{k+1}j_{k+1}}\ldots A_{i_{l-1}j_{l-1}})d(A_{i_kj_k} A_{i_lj_l})\ldots\end{eqnarray*}
By applying Itô's integration by part formula to quantities of the kind $\sum_s \qU_{si}^*\qU_{sj}$ or $\sum_s \qU_{is}\qU_{js}^*$, we find that the differential of this quantities needs to be zero. Taking into account the initial condition, we find that both quantities equal $\delta_{ij}$. In other words, the process $\qU(t)$ can be seen as a stochastic process on $U\langle n\rangle$ via the map:
\begin{eqnarray*}
\mathfrak j&:&U_n\rightarrow A\\
&&u_{ij}\mapsto \qU_{ij}
\end{eqnarray*}
\\
We will consider quantities of the type $\phi(\qU_{i_{11}j_{11}}^{\epsilon_{11}}\ldots \qU_{i_{k_11}j_{k_11}}^{\epsilon_{k_11}})\ldots\phi(\qU_{i_{1r}j_{1r}}^{\epsilon_{1r}}\ldots \qU_{i_{k_1r}j_{k_1r}}^{\epsilon_{k_1r}})$. We will need a convenient notation to denote such functions in a way that is consistent with the transformation they will undergo. Thus, we will denote for instance by:
$$T=\begin{array}{|c|c|c|c|c||c|c|c||c||c|c|}\hline \multicolumn{5}{|c||}{1}&\multicolumn{3}{c||}{2}&3&4\\\hline \mathfrak U&\mathfrak U&\mathfrak U&\mathfrak U&\mathfrak U&\mathfrak U&\ldots&\mathfrak U&\mathfrak U&\mathfrak U\\(1,2)&(2,3)&(5,2)&(1,1)&(2,1)&(1,1)&\ldots&(4,5)&(1,1)&(5,5)\\1&0&0&1&0&0&\ldots&1&1&0\\\hline\end{array}$$
the function:
$$\phi(\mathfrak U^*_{1,2}\mathfrak U_{2,3}\mathfrak U_{5,2}\mathfrak U^*_{1,1}\mathfrak U_{2,1})\phi(\mathfrak U_{1,1}\ldots\mathfrak U^*{4,5})\phi(\mathfrak U_{1,1}^*)\phi(\mathfrak U_{5,5})$$
In other words, the first line gives the rank of the $\phi$ we are in, the second line (which is frankly not necessary) reminds us of the type of matricial blocks we consider (here blocks of $\qU$), the third line gives the coordinates of the blocks, and the last lines tells us whether or not there is a $*$.\\
We will call this kind of functions, \''table-functions\'', because of the way they are denoted.\\
Then, we can write:
$$T'=\hat{\mathcal D} T+\hat{\mathcal U} T$$
where $\hat{\mathcal D}$ is the operator corresponding to the term with one differential, $\hat{\mathcal U}$ is the term corresponding to the quadratic variation inside a same $\phi$. Let us remark that the quadratic variation can only affect one and the same trace.  Let us examine these operators in detail. Because we take the expectation $\phi$ it must be remembered that the (quantum) martingale part vanishes and we therefore do not explicit it.
\subsection{The operator $\hat{\mathcal D}$:}it is easy to see that it behaves in the following way:

\begin{figure}[H]
\centering
\begin{tikzpicture}
\draw (0,0) node[above]{$\begin{array}{|c||c|c|c||c|}\hline\ldots&\multicolumn{3}{c||}{\kappa}&\ldots\\\hline\ldots&\ldots&\mathfrak U&\ldots&\ldots\\\hline\ldots&\ldots&(i_{\mu \kappa},j_{\mu\kappa})&\ldots&\ldots\\\hline\ldots&\ldots&0&\ldots&\ldots\\\hline\end{array}$};
\draw [-stealth](0,0) -- (0,-2);
\draw (0,-2) node[below]{$\frac12n\sum_{\kappa=1}^r\sum_{1\leq\mu\leq k_\kappa}\begin{array}{|c||c|c|c||c|}\hline\ldots&\multicolumn{3}{c||}{\kappa}&\ldots\\\hline\ldots&\ldots&\mathfrak U&\ldots&\ldots\\\hline\ldots&\ldots&(i_{\mu \kappa},j_{\mu\kappa})&\ldots&\ldots\\\hline\ldots&\ldots&0&\ldots&\ldots\\\hline\end{array}$};
\end{tikzpicture}
\end{figure}
and
\begin{figure}[H]
\centering
\begin{tikzpicture}
\draw (0,0) node[above]{$\begin{array}{|c||c|c|c||c|}\hline\ldots&\multicolumn{3}{c||}{\kappa}&\ldots\\\hline\ldots&\ldots&\mathfrak U&\ldots&\ldots\\\hline\ldots&\ldots&(i_{\mu \kappa},j_{\mu\kappa})&\ldots&\ldots\\\hline\ldots&\ldots&1&\ldots&\ldots\\\hline\end{array}$};
\draw [-stealth](0,0) -- (0,-2);
\draw (0,-2) node[below]{$\frac12n\sum_{\kappa=1}^r\sum_{1\leq\mu\leq k_\kappa}\begin{array}{|c||c|c|c||c|}\hline\ldots&\multicolumn{3}{c||}{\kappa}&\ldots\\\hline\ldots&\ldots&\mathfrak U&\ldots&\ldots\\\hline\ldots&\ldots&(i_{\mu \kappa},j_{\mu\kappa})&\ldots&\ldots\\\hline\ldots&\ldots&1&\ldots&\ldots\\\hline\end{array}$};
\end{tikzpicture}
\end{figure}
In other words, we obtain $\frac12\left(k_1+\ldots+k_\kappa\right)T$.
\subsection{The operator $\hat{\mathcal U}$:} we have the following behavior.
\begin{figure}[H]
\centering
\begin{tikzpicture}
\draw (0,0) node[above]{$\begin{array}{|c||c|c|c|c|c||c|}\hline\ldots&\multicolumn{5}{c||}{\kappa}&\ldots\\\hline\ldots&\ldots&\mathfrak U&\ldots&\mathfrak U&\ldots&\ldots\\\hline\ldots&\ldots&(i_{\mu_1 \kappa},j_{\mu_1\kappa})&\ldots&(i_{\mu_2\kappa},j_{\mu_2\kappa})&\ldots&\ldots\\\hline\ldots&\ldots&0&\ldots&0&\ldots&\ldots\\\hline\end{array}$};
\draw [-stealth](0,0) -- (0,-2);
\draw (0,-2) node[below]{$-\sum_{\kappa=1}^r\sum_{1\leq\mu_1<\mu_2\leq k_\kappa}\begin{array}{|c||c|c|c|c|c||c||c|c|c|}\hline\ldots&\multicolumn{5}{c||}{\kappa}&\ldots&\multicolumn{3}{c|}{r+1}\\\hline\ldots&\ldots&&\mathfrak U&&\ldots&\ldots&&\ldots&\mathfrak U\\\hline\ldots&\ldots&\alpha_{\mu_1-1,\kappa}&(i_{\mu_1 \kappa},j_{\mu_2\kappa})&\alpha_{\mu_2+1,\kappa}&\ldots&\ldots&\alpha_{\mu_1+1,\kappa}&\ldots&(i_{\mu_2\kappa},j_{\mu_1,\kappa})\\\hline\ldots&\ldots&&0&&\ldots&\ldots&&\ldots&0\\\hline\end{array}$};

\end{tikzpicture}
\end{figure}
and
\begin{figure}[H]
\centering
\begin{tikzpicture}
\draw (0,0) node[above]{$\begin{array}{|c||c|c|c|c|c||c|}\hline\ldots&\multicolumn{5}{c||}{\kappa}&\ldots\\\hline\ldots&\ldots&\mathfrak U&\ldots&\mathfrak U&\ldots&\ldots\\\hline\ldots&\ldots&(i_{\mu_1 \kappa},j_{\mu_1\kappa})&\ldots&(i_{\mu_2\kappa},j_{\mu_2\kappa})&\ldots&\ldots\\\hline\ldots&\ldots&1&\ldots&1&\ldots&\ldots\\\hline\end{array}$};
\draw [-stealth](0,0) -- (0,-2);
\draw (0,-2) node[below]{$-\sum_{\kappa=1}^r\sum_{1\leq\mu_1<\mu_2\leq k_\kappa}\begin{array}{|c||c|c|c|c|c||c||c|c|c|}\hline\ldots&\multicolumn{5}{c||}{\kappa}&\ldots&\multicolumn{3}{c|}{r+1}\\\hline\ldots&\ldots&&\mathfrak U&&\ldots&\ldots&&\ldots&\mathfrak U\\\hline\ldots&\ldots&\alpha_{\mu_1-1,\kappa}&(i_{\mu_2 \kappa},j_{\mu_1\kappa})&\alpha_{\mu_2+1,\kappa}&\ldots&\ldots&\alpha_{\mu_1+1,\kappa}&\ldots&(i_{\mu_1\kappa},j_{\mu_2,\kappa})\\\hline\ldots&\ldots&&1&&\ldots&\ldots&&\ldots&1\\\hline\end{array}$};
\end{tikzpicture}
\end{figure}
and
\begin{figure}[H]
\centering
\begin{tikzpicture}
\draw (0,0) node[above]{$\begin{array}{|c||c|c|c|c|c||c|}\hline\ldots&\multicolumn{5}{c||}{\kappa}&\ldots\\\hline\ldots&\ldots&\mathfrak U&\ldots&\mathfrak U&\ldots&\ldots\\\hline\ldots&\ldots&(i_{\mu_1 \kappa},j_{\mu_1\kappa})&\ldots&(i_{\mu_2\kappa},j_{\mu_2\kappa})&\ldots&\ldots\\\hline\ldots&\ldots&0&\ldots&1&\ldots&\ldots\\\hline\end{array}$};
\draw [-stealth](0,0) -- (0,-2);
\draw (-1,-2) node[below]{$-\delta_{j_{\mu_1\kappa}j_{\mu_2\kappa}}\sum_s\sum_\kappa\sum_{\mu_1<\mu_2}\begin{array}{|c||c|c|c|c|c|c||c||c|c|c|}\hline\ldots&\multicolumn{6}{c||}{\kappa}&\ldots&\multicolumn{3}{c|}{r+1}\\\hline\ldots&\ldots&\qU&\mathfrak U&\qU&\qU&\ldots&\ldots&&\ldots&\mathfrak U\\\hline\ldots&\ldots&\alpha_{\mu_1-1,\kappa}&(i_{\mu_1 \kappa},s)&(i_{\mu_2\kappa},s)&\alpha_{\mu_2+1,\kappa}&\ldots&\ldots&\alpha_{\mu_1+1,\kappa}&\ldots&(i_{\mu_2\kappa},j_{\mu_1,\kappa})\\\hline\ldots&\ldots&&0&1&&\ldots&\ldots&&\ldots&\\\hline\end{array}$};
\end{tikzpicture}
\end{figure}
and, finally,
\begin{figure}[H]
\centering
\begin{tikzpicture}
\draw (0,0) node[above]{$\begin{array}{|c||c|c|c|c|c||c|}\hline\ldots&\multicolumn{5}{c||}{\kappa}&\ldots\\\hline\ldots&\ldots&\mathfrak U&\ldots&\mathfrak U&\ldots&\ldots\\\hline\ldots&\ldots&(i_{\mu_1 \kappa},j_{\mu_1\kappa})&\ldots&(i_{\mu_2\kappa},j_{\mu_2\kappa})&\ldots&\ldots\\\hline\ldots&\ldots&1&\ldots&0&\ldots&\ldots\\\hline\end{array}$};
\draw [-stealth](0,0) -- (0,-2);
\draw (-1,-2) node[below]{$-\delta_{j_{\mu_1\kappa}j_{\mu_2\kappa}}\sum_s\sum_\kappa\sum_{\mu_1<\mu_2}\begin{array}{|c||c|c|c|c||c||c|c|c|}\hline\ldots&\multicolumn{4}{c||}{\kappa}&\ldots&\multicolumn{3}{c|}{r+1}\\\hline\ldots&\ldots&\qU&\mathfrak U&\ldots&\ldots&\qU&\ldots&\mathfrak U\\\hline\ldots&\ldots&\alpha_{\mu_1-1,\kappa}&\alpha_{\mu_2+1,\kappa}&\ldots&\ldots&(i_{\mu_1\kappa},s)&\ldots&(i_{\mu_2\kappa},s)\\\hline\ldots&\ldots&&&\ldots&\ldots&1&\ldots&0\\\hline\end{array}$};
\end{tikzpicture}
\end{figure}
We notice that we obtain a very similar result to what was obtained at the end of section 3 in \cite{Ulrich2015}, with the difference that the summation over $s$ is made on the second index here, instead of the first one. This is due to our having taken the left version of the stochastic equation instead of the right one.

\section{Convergence of the Brownian motion on the orthogonal group}
\subsection{The theorem}
In the sequel of this section, we will always denote by $n$ and $m$ two (natural) integers. The first one, $n$, will be fixed, and the second one, $m$ will tend to infinity, so that we can study the blockwise structure of the limit.\\
We want to study the high-dimensional limit of the Brownian motion on $O(nm)$, so let us first give its stochastic differential equation. We follow \cite{Dahlqvist2017}. Let us note that we actually follow the indication to build the Brownian motion on $O(nm)$, but because this Lie group is not connected and we begin the process at $I_{nm}$, it is the Brownian motion on $SO(nm)$ that we get.\\
Let $(X(t))_{t\geq0}$ be a matricial stochastic process build in the following way:
\begin{itemize}
    \item The family $(X_{k,l}(t))_{1\leq k\leq l\leq nm}$ is a family of independent variables
    \item The family $(X_{k,l}(t))_{1\leq k< l\leq nm}$ is a family of i.i.d. random variables following a Gaussian law $\mathcal N(0,t/m)$.
    \item $X_{kk}(t)=0$ for each $k$
    \item The matrix $X(t)$ is anti-symmetric.
\end{itemize}
It is then clear that $(X(t))_t$ is a Brownian motion on the space of anti-symmetric matrices, which is the Lie algebra of $O(nm)$.\\
We will need the quadratic variation of the $X_{j,k}$, so let us compute them right away:

$$d[X_{ab},X_{cd}]=\left\{\begin{array}{cc}dt/m&\text{if }a=c\text{ and }b=d\text{ and }a\neq b\\-dt/m&\text{ if }a=d\text{ and }b=c\text{ and }a\neq b\\0&\text{ else}\end{array}\right.$$
The quadratic variation matrix mentioned in \cite{Dahlqvist2017} needs to be understood as the matrix whose coefficients are the quadratic variations of $ dX_t.dX_t$, seen as a product of matrices. It is therefore: $$\langle dX_t.dX_t\rangle_{j,k}=\sum_{k=1}^{nm}d[X_{jl},X_{lk}]=\delta_{jk}\sum_{k=1,k\neq j}^{nm}\frac{dt}{m}=\delta_{jk}dt\frac{nm-1}m$$
This matrix thus is $\frac{nm-1}{m}dt I_{nm}$, and so the stochastic differential equation for 

the Brownian motion on $SO(mn)$ is:
$$dO_t=O_tdX_t+\frac12(n-\frac1m)O_tdt$$
with initial condition $O_0=I_{nm}$.\\
Let us note that we obtain a consistent result with \cite{Levy2011}. We only need to remember that we did not take the same renormalization\footnote{Let us remark also that this means that in \cite{Ulrich2015} I forgot to take the renormalization into account in the equation. Indeed, the last term should have been $-ndtU_tdt$ instead of just $-U_tdt$. This does not change the result, as the whole proof can easily be adapted. Of course, \cite{CU2016}, which obtains the same result by an other method, also confirms the result.} than Lévy, which explains why the factor in front of the $dt$ term is not the same. \\
We can also remark that the most natural framework to express a Brownian motion on a Lie group is the Stratonovich integral. Indeed, the formula in \cite{Dahlqvist2017} for such a Brownian motion, $dG_t=G_tdK_t+\frac{C_{\mathfrak g}}{2}dt$, together with the Itô-Stratonovich conversion formula ($d\circ M_t=dM_t+\frac12d[M_t,M_t]$ for any process $M_t$ for which this has a sense), implies that this Brownian motion can be described by the Stratonovich equation $dG_t=G_t d\circ K_t$. Nevertheless, such an equation is less tractable for the computations. Indeed, the advantage of Itô's formula is the fact that an element of the kind $G_t dK_t$ is always a martingale, and thus vanishes when we take the expectation. This is why we favor Itô's framework.\\
Finally, let us remark that if we compute $d(O_tO_t^*)$ and $d(O_t^*O_t)$ through Itô's integration by part formula we find in both cases zero. This is consistent with the fact that $O_t$ remains on $O(nm)$.
\\
We will consider that all the (classical) random variables considered thus far are defined on an unique probability space $(\Omega,\mathcal F,\mathbb P)$. We will also denote by $L^{\infty_-}(\Omega)$ the space of random variables defined on $\Omega$ admitting moments of every order.\\

Using the dual semigroups introduced in the previous section, we can define:
$$\begin{array}{cc}j_m(t):&A\rightarrow L^{\infty_-}(\Omega)\\&u_{ij}\mapsto [O(t)]_{ij}\end{array}$$
where $O(t)$ is the orthogonal Brownian motion on $O(nm)$ and $[O(t)]_{i,j}$ designates the $(i,j)$ block of size $m\times m$ of this process.

We will prove the following theorem in this section:
\begin{thm}[Convergence of the orthogonal Brownian motion]\label{CVortho}
When $m$ tends to the infinite, the processes $j_m$ tend in $*$-moments towards the process $\mathfrak j$.
\end{thm}

\subsection{Notations}
Before going over to the proof, though, we need to introduce some useful notations that will serve in the proof.\\
We introduce $[\mathcal J]=\{1,\ldots,n\}^2$, $\mathcal J=\{1,\ldots,mn\}^2$. For a matrix $M$ of size $mn\times mn$, we will denote by $[M]_{kl}$ the $(k,l)$-block of size $m\times m$ (for $(k,l)\in [\mathcal J]$), and by $M_{kl}$ the coefficient $(k,l)$ of the matrix (for $(k,l)\in \mathcal J$). 
Whenever we write $Tr(A)$, we intend to mean the trace of matrix $A$, but when we write $tr(A)$, we actually intend $tr(A)=\frac1m Tr(A)$.\\
We also emphasize the difference between $[O^*]_{i,j}$ and $[O]_{ij}^*$: we have the relationship: $([O]^*_{i,j})_{a,b}=([O]_{i,j})_{b,a}$, where $([M]_{i,j})_{a,b}$ denotes the $(a,b)$ coefficient of the block $(i,j)$ of matrix $M$.\\
As in \cite{Ulrich2015}, we will be interested in the differential equations satisfied by some functions. To denote the functions, we will use tables, in the following fashion:
$$T=\begin{array}{|c|c|c|c|c||c|c|c||c||c|c|}\hline \multicolumn{5}{|c||}{1}&\multicolumn{3}{c||}{2}&3&4\\\hline O&O&O&O&O&O&\ldots&O&O&O\\(1,2)&(2,3)&(5,2)&(1,1)&(2,1)&(1,1)&\ldots&(4,5)&(1,1)&(5,5)\\1&0&0&1&0&0&\ldots&1&1&0\\\hline\end{array}$$
Then, $T$ denotes the function $\mathbb{E}[tr([O]^*_{1,2}[O]_{2,3}[O]_{5,2}[O]^*_{1,1}[O]_{2,1})tr([O]_{1,1}\ldots[O]^*_{4,5})tr([O]_{1,1}^*)tr([O]_{5,5})]$. We emphasize the fact that we use the table as a name for the function. Thus, if $T$ is this table, we will be able to write things like $T^\prime$, etc.

\subsection{The differential equation for the classical process}
We will need to apply Itô formula later, so we can already compute some useful quantities. 
We observe that, due to what we know about the quadratic variation of the $X_{j,k}$, we have:
\begin{eqnarray*}d[O_{ab},O_{cd}]&=&\sum_{s=1}^{mn}\sum_{q=1}^{mn}O_{as}O_{cq}d[X_{sb},X_{qd}]\\&=&\underbrace{\delta_{bd}\sum_{s=1,s\neq b}^{mn}O_{as}O_{cs}dt/m}_{(1)}-\underbrace{(1-\delta_{bd})O_{ad}O_{cb}dt/m}_{(2)}\end{eqnarray*}
We see that this quadratic variation can be decomposed into two parts, $(1)$ and $(2)$. \\
When considering a function $T$, we are actually considering sums of terms of the type:
\begin{eqnarray*}(1)M^{\epsilon_{11}}_{(j_{11}-1)m+s_{11},(l_{11}-1)m+s_{21}}&&M^{\epsilon_{21}}_{(j_{21}-1)m+s_{21},(l_{21}-1)m+s_{31}}\ldots M^{\epsilon_{k_11}}_{(j_{k_11}-1)m+s_{k_11},(l_{k_11}-1)m+s_{11}}M^{\epsilon_{12}}_{(j_{12}-1)m+s_{12},(l_{12}-1)m+s_{22}}\times\\&&\ldots  M^{\epsilon_{k_22}}_{(j_{k_22}-1)m+s_{k_22},(l_{k_22}-1)m+s_{12}}\ldots\end{eqnarray*}
In particular, we observe how the $s$'s work, due to the fact that we are considering traces. 
\\
To find a differential equation verified by $T$, we need to apply Itô's formula to calculate the differential of an expression of the type $(1)$. When computing this formula, we are not interested in the martingal parts, because they vanish once we take the expectation. Thus we have something of the kind:
\begin{eqnarray*}d(\ldots)&=&\underbrace{\sum_{\kappa=1}^r\sum_{\mu=1}^{k_\kappa}\ldots d(O^{\epsilon_{\mu\kappa}}_{(j_{\mu\kappa}-1)m+s_{\mu\kappa},(l_{\mu\kappa}m-1)+s_{\mu+1,\kappa}})\ldots}_{\mathcal{D}}\\&+&\underbrace{\sum_{\kappa=1}^r\sum_{1\leq \mu_1<\mu_2\leq k_\kappa}\ldots d[O^{\epsilon_{\mu_1\kappa}}_{(j_{\mu_1\kappa}-1)m+s_{\mu_1\kappa},(l_{\mu_1\kappa}-1)m+s_{\mu_1+1,\kappa}},O^{\epsilon_{\mu_2\kappa}}_{(j_{\mu_2\kappa}-1)m+s_{\mu_2\kappa},(l_{\mu_2\kappa}-1)m+s_{\mu_2+1,\kappa}}]\ldots}_{\mathcal U}\\&+&\underbrace{\sum_{1\leq\kappa_1<\kappa_2\leq r}\sum_{\mu_1=1}^{k_{\kappa_1}}\sum_{\mu_2=1}^{k_{\kappa_2}}\ldots d[O^{\epsilon_{\mu_1\kappa_1}}_{(j_{\mu_1\kappa_1}-1)m+s_{\mu_1\kappa_1},(l_{\mu_1\kappa_1}-1)m+s_{\mu_1+1,\kappa_1}},O^{\epsilon_{\mu_2\kappa_2}}_{(j_{\mu_2\kappa_2}-1)m+s_{\mu_2\kappa_2},(l_{\mu_2\kappa_2}-1)m+s_{\mu_2+1,\kappa_2}}]\ldots}_{\mathcal B}\end{eqnarray*}
We see that three different operators naturally appear, namely $\mathcal D$, $\mathcal U$ and $\mathcal B$. When remembering that $d[O_{ab},O_{cd}]$ has two terms, we come up with five operators: $\mathcal D$, $\mathcal U^{(1)}$, $\mathcal U^{(2)}$, $\mathcal B^{(1)}$ and $\mathcal B^{(2)}$. We will examine in turn their actions in the following. But let us first make two remarks. First, we need to keep in mind that, while evaluating the quantities, the martingale terms vanish because of the expectation. Second, we will see that we will arrive at a system of differential equations. As the various operators create or erase traces but do not change the overall number of blocks of $O$ intervening in the function, it must be noted that the system characterizing a specific function $T$ makes use of only a finite number of other functions obtained from $T$ by combinatorial means. More precisely, let us call the \emph{order} of $T$ the number of matricial blocks coming into play in $T$. This means that the order $\eta$ is $\eta=k_1+\ldots+k_r$. From Itô's formula combined with the stochastic equations, it is clear that no additional matricial block can appear when you derivate $T$. Neither can a block disappear. So, we see that all different functions intervening in the expression of the derivative $T'$ have all the same order $\eta$. This in turn means that the equation $T'$ is part of a system of differential equation that has at most $n^{2\eta}2^\eta\sum_{z=1}^\eta \binom{\eta-1}{z-1}$ functions and $n^{2\eta}2^\eta\sum_{z=1}^\eta \binom{\eta-1}{z-1}$ equations. This is computed considering that the indices of a matricial block can vary, as well as the fact of having a $*$ or not, and that the number of traces and where they begin and where they end can be modified. Though this quantity may be quite large, it is nonetheless finite.
 Thus, if we denote by $B$ the $n^{2\eta}2^\eta\sum_{z=1}^\eta \binom{\eta-1}{z-1}\times n^{2\eta}2^\eta\sum_{z=1}^\eta \binom{\eta-1}{z-1}$ matrix of coefficients for this system of differential equation, the expression of $T_m(t)$ will be one of the lines of the column-vector $e^{Bt}v_0$, where $v_0$ is a vector containing the initial conditions.\\
 Now, we will see that the equation for $T'$ is of the type $T'=C+D$, where $C$ contains only constant coefficients in front of the functions and $D$ contains coefficients in front of the functions that are in $O(\frac1m)$. Therefore, the matrix $B$ can be decomposed into $B=B_0+O(\frac1m)$, where $B_0$ contains only constants. Hence, the expression of $T_m(t)$, when $m$ tends to infinity, will tend towards $e^{B_0t}v_0$. In other words, $T_m(t)$ will tend pointwise, in $m$, towards a function $\tilde T$ that verifies the equation $\tilde T'=C$.\\
 In practice, this means that while studying the four operators, we will be interested only in those parts of the operators that are not in $O(\frac1m)$. We will denote by $\mathcal D_\infty$, $\mathcal U^{(1)}_\infty$, and so forth, the part of the respective operator that is constant in $m$.

\subsubsection{The operator $\mathcal D$:}The operator is characterized by the following behaviors:
\begin{figure}[H]
\centering
\begin{tikzpicture}
\draw (0,0) node[above]{$\begin{array}{|c||c|c|c||c|}\hline\ldots&\multicolumn{3}{c||}{\kappa}&\ldots\\\hline\ldots&\ldots&O&\ldots&\ldots\\\hline\ldots&\ldots&(i_{\mu \kappa},j_{\mu\kappa})&\ldots&\ldots\\\hline\ldots&\ldots&0&\ldots&\ldots\\\hline\end{array}$};
\draw [-stealth](0,0) -- (0,-2);
\draw (0,-2) node[below]{$\frac12(n-\frac1m)\sum_{\kappa=1}^r\sum_{1\leq\mu\leq k_\kappa}\begin{array}{|c||c|c|c||c|}\hline\ldots&\multicolumn{3}{c||}{\kappa}&\ldots\\\hline\ldots&\ldots&O&\ldots&\ldots\\\hline\ldots&\ldots&(i_{\mu \kappa},j_{\mu\kappa})&\ldots&\ldots\\\hline\ldots&\ldots&0&\ldots&\ldots\\\hline\end{array}$};
\end{tikzpicture}
\end{figure}
and
\begin{figure}[H]
\centering
\begin{tikzpicture}
\draw (0,0) node[above]{$\begin{array}{|c||c|c|c||c|}\hline\ldots&\multicolumn{3}{c||}{\kappa}&\ldots\\\hline\ldots&\ldots&O&\ldots&\ldots\\\hline\ldots&\ldots&(i_{\mu \kappa},j_{\mu\kappa})&\ldots&\ldots\\\hline\ldots&\ldots&1&\ldots&\ldots\\\hline\end{array}$};
\draw [-stealth](0,0) -- (0,-2);
\draw (0,-2) node[below]{$\frac12(n-\frac1m)\sum_{\kappa=1}^r\sum_{1\leq\mu\leq k_\kappa}\begin{array}{|c||c|c|c||c|}\hline\ldots&\multicolumn{3}{c||}{\kappa}&\ldots\\\hline\ldots&\ldots&O&\ldots&\ldots\\\hline\ldots&\ldots&(i_{\mu \kappa},j_{\mu\kappa})&\ldots&\ldots\\\hline\ldots&\ldots&1&\ldots&\ldots\\\hline\end{array}$};
\end{tikzpicture}
\end{figure}
In other words, we obtain always $\frac12(n-\frac 1m)(k_1+\ldots+k_r)T$. We see that this does tend, when $m$ goes to infinity, towards $\frac n2(k_1+\ldots+k_r)T$, which corresponds to $\hat{\mathcal D}$.

\subsubsection{The operator $\mathcal U^{(2)}$:} We proceed as in the previous case, but here we need to distinguish four cases, according to whether or not $\epsilon_{\mu_1\kappa}$ and $\epsilon_{\mu_2\kappa}$ are zero or one. We notice that we have almost the same behavior as for the Brownian motion on the unitary group, as described in \cite{Ulrich2015}. The only real difference comes from the $1-\delta_{bd}$ factor in the quadratic variation. Let us examine in some details the case where $\epsilon_{\mu_1\kappa}=\epsilon_{\mu_2\kappa}=0$, ie when we have no $*$ on the two blocks involved in the quadratic variation. We will then be able to treat the other cases faster. Because of the $1-\delta_{bd}$ factor, if $j_{\mu_1}\neq j_{\mu_2}$, then this factor vanishes and we have only:

\begin{figure}[H]
\centering
\begin{tikzpicture}
\draw (0,0) node[above]{$\begin{array}{|c||c|c|c|c|c||c|}\hline\ldots&\multicolumn{5}{c||}{\kappa}&\ldots\\\hline\ldots&\ldots&O&\ldots&O&\ldots&\ldots\\\hline\ldots&\ldots&(i_{\mu_1 \kappa},j_{\mu_1\kappa})&\ldots&(i_{\mu_2\kappa},j_{\mu_2\kappa})&\ldots&\ldots\\\hline\ldots&\ldots&0&\ldots&0&\ldots&\ldots\\\hline\end{array}$};
\draw [-stealth](0,0) -- (0,-2);
\draw (0,-2) node[below]{$-\sum_{\kappa=1}^r\sum_{1\leq\mu_1<\mu_2\leq k_\kappa}\begin{array}{|c||c|c|c|c|c||c||c|c|c|}\hline\ldots&\multicolumn{5}{c||}{\kappa}&\ldots&\multicolumn{3}{c|}{r+1}\\\hline\ldots&\ldots&&O&&\ldots&\ldots&&\ldots&O\\\hline\ldots&\ldots&\alpha_{\mu_1-1,\kappa}&(i_{\mu_1 \kappa},j_{\mu_2\kappa})&\alpha_{\mu_2+1,\kappa}&\ldots&\ldots&\alpha_{\mu_1+1,\kappa}&\ldots&(i_{\mu_2\kappa},j_{\mu_1,\kappa})\\\hline\ldots&\ldots&&0&&\ldots&\ldots&&\ldots&0\\\hline\end{array}$};

\end{tikzpicture}
\end{figure}
This is completely similar to what was obtained in \cite{Ulrich2015}. In particular, we observe that even though we had a factor $1/m$, because a new trace appears and because we use the trace normalized with $1/m$, we have only constant factors here. Now, what happens if $j_{\mu_1}= j_{\mu_2}$? We obtain almost the same result, with the exception that all the cases where $s_{\mu_1+1,\kappa}=s_{\mu_2+1,\kappa}$ vanish. The difference between the expected
$$-\sum_{\kappa=1}^r\sum_{1\leq\mu_1<\mu_2\leq k_\kappa}\begin{array}{|c||c|c|c|c|c||c||c|c|c|}\hline\ldots&\multicolumn{5}{c||}{\kappa}&\ldots&\multicolumn{3}{c|}{r+1}\\\hline\ldots&\ldots&&O&&\ldots&\ldots&&\ldots&O\\\hline\ldots&\ldots&\alpha_{\mu_1-1,\kappa}&(i_{\mu_1 \kappa},j_{\mu_2\kappa})&\alpha_{\mu_2+1,\kappa}&\ldots&\ldots&\alpha_{\mu_1+1,\kappa}&\ldots&(i_{\mu_2\kappa},j_{\mu_1,\kappa})\\\hline\ldots&\ldots&&0&&\ldots&\ldots&&\ldots&0\\\hline\end{array}$$
and what we really obtain due to the vanishing of these terms, is:
$$\frac{dt}{m}\sum_{\dots}\sum_{k=1}^m\ldots ([O]_{i_{\mu_1\kappa}j_{\mu_2\kappa}})_{s_{\mu_1\kappa}k}\ldots([O]_{i_{\mu_2\kappa}j_{\mu_1\kappa}})_{s_{\mu_2\kappa}k}$$
This is the same as:
\begin{eqnarray*}\frac{dt}{m}\sum_{\dots}\sum_{k=1}^m\ldots ([O]_{i_{\mu_1\kappa}j_{\mu_2\kappa}})_{s_{\mu_1\kappa}k}&&([O]^*_{i_{\mu_2\kappa}j_{\mu_1\kappa}})_{ks_{\mu_2\kappa}}([O]^{1-\epsilon_{\mu_2-1,\kappa}}_{i_{\mu_2-1,\kappa}j_{\mu_2-1,\kappa}})_{s_{\mu_2\kappa}s_{\mu_2-1,\kappa}}\ldots([O]^{1-\epsilon_{\mu_1+1,\kappa}}_{i_{\mu_1+1,\kappa}j_{\mu_1+1,\kappa}})_{s_{\mu_1+2\kappa}s_{\mu_1+1,\kappa}}\\
&&\times([O]_{i_{\mu_2+1\kappa}j_{\mu_2+1\kappa}}^{\epsilon_{\mu_2+1\kappa}})_{s_{mu_2+1},s_{\mu_2+2}}\ldots\end{eqnarray*}
This is due to the fact that $s_{\mu_1+1}=s_{\mu_2+1}$. So, at the price of interverting $*$ and no-$*$ for a part of the trace, we were able to write this difference as having the exact same number of traces as the original $T$. But because there is the factor $\frac1m$ in front, this difference vanishes in the $m\rightarrow\infty$ limit. So, we have in all case, that the operator $\mathcal U^{(2)}_\infty$ acts as such:

\begin{figure}[H]
\centering
\begin{tikzpicture}
\draw (0,0) node[above]{$\begin{array}{|c||c|c|c|c|c||c|}\hline\ldots&\multicolumn{5}{c||}{\kappa}&\ldots\\\hline\ldots&\ldots&O&\ldots&O&\ldots&\ldots\\\hline\ldots&\ldots&(i_{\mu_1 \kappa},j_{\mu_1\kappa})&\ldots&(i_{\mu_2\kappa},j_{\mu_2\kappa})&\ldots&\ldots\\\hline\ldots&\ldots&0&\ldots&0&\ldots&\ldots\\\hline\end{array}$};
\draw [-stealth](0,0) -- (0,-2);
\draw (0,-2) node[below]{$-\sum_{\kappa=1}^r\sum_{1\leq\mu_1<\mu_2\leq k_\kappa}\begin{array}{|c||c|c|c|c|c||c||c|c|c|}\hline\ldots&\multicolumn{5}{c||}{\kappa}&\ldots&\multicolumn{3}{c|}{r+1}\\\hline\ldots&\ldots&&O&&\ldots&\ldots&&\ldots&O\\\hline\ldots&\ldots&\alpha_{\mu_1-1,\kappa}&(i_{\mu_1 \kappa},j_{\mu_2\kappa})&\alpha_{\mu_2+1,\kappa}&\ldots&\ldots&\alpha_{\mu_1+1,\kappa}&\ldots&(i_{\mu_2\kappa},j_{\mu_1,\kappa})\\\hline\ldots&\ldots&&0&&\ldots&\ldots&&\ldots&0\\\hline\end{array}$};
\end{tikzpicture}
\end{figure}

We also need to treat the case when one or two of the $O$ appearing in the quadratic variation are decorated with a $*$. Let us treat it in some details when both have the $*$. The other cases will be done in a similar fashion. We first observe that $([O]_{jl}^*)_{s_1s_2}=([O]_{jl})_{s_2s_1}$, the star being only a transposition. So, we are considering products of the kind:
$$([O]_{i_1j_1}^{\epsilon_1})_{s_1s_2}\ldots ([O]_{i_{\mu_1}j_{\mu_1}}^*)_{s_{\mu_1}s_{\mu_1+1}}\ldots ([O]_{i_{\mu_2}j_{\mu_2}}^*)_{s_{\mu_2}s_{\mu_2+1}}\ldots ([O]_{i_{k_\kappa}j_{k_\kappa}}^{\epsilon_{k_\kappa}})_{s_{k_\kappa}s_1}$$
This is then the same as:
$$([O]_{i_1j_1}^{\epsilon_1})_{s_1s_2}\ldots ([O]_{i_{\mu_1}j_{\mu_1}})_{s_{\mu_1+1}s_{\mu_1}}\ldots ([O]_{i_{\mu_2}j_{\mu_2}})_{s_{\mu_2+1}s_{\mu_2}}\ldots ([O]_{i_{k_\kappa}j_{k_\kappa}}^{\epsilon_{k_\kappa}})_{s_{k_\kappa}s_1}$$
When we apply the operator, because of how the quadratic variation behaves, when get:
$$([O]_{i_1j_1}^{\epsilon_1})_{s_1s_2}\ldots ([O]_{i_{\mu_1}j_{\mu_2}})_{s_{\mu_1+1}s_{\mu_2}}\ldots ([O]_{i_{\mu_2}j_{\mu_1}})_{s_{\mu_2+1}s_{\mu_1}}\ldots ([O]_{i_{k_\kappa}j_{k_\kappa}}^{\epsilon_{k_\kappa}})_{s_{k_\kappa}s_1}$$
which is the same as
$$([O]_{i_1j_1}^{\epsilon_1})_{s_1s_2}\ldots ([O]_{i_{\mu_1}j_{\mu_2}}^*)_{s_{\mu_2}s_{\mu_1+1}}\ldots ([O]^*_{i_{\mu_2}j_{\mu_1}})_{s_{\mu_1}s_{\mu_2+1}}\ldots ([O]_{i_{k_\kappa}j_{k_\kappa}}^{\epsilon_{k_\kappa}})_{s_{k_\kappa}s_1}$$
If we wanted to treat it comprehensively, we would also need to distinguish whether or not $j_{\mu_1\kappa}=j_{\mu_2\kappa}$. But when this is the case, the difference between what we obtain and what we really get because of vanishing terms can be treated as before and is again a $O(\frac1m)$. So, the operator $\mathcal U^{(2)}_\infty$ acts on the tables in the following way:
\begin{figure}[H]
\centering
\begin{tikzpicture}
\draw (0,0) node[above]{$\begin{array}{|c||c|c|c|c|c||c|}\hline\ldots&\multicolumn{5}{c||}{\kappa}&\ldots\\\hline\ldots&\ldots&O&\ldots&O&\ldots&\ldots\\\hline\ldots&\ldots&(i_{\mu_1 \kappa},j_{\mu_1\kappa})&\ldots&(i_{\mu_2\kappa},j_{\mu_2\kappa})&\ldots&\ldots\\\hline\ldots&\ldots&1&\ldots&1&\ldots&\ldots\\\hline\end{array}$};
\draw [-stealth](0,0) -- (0,-2);
\draw (0,-2) node[below]{$-\sum_{\kappa=1}^r\sum_{1\leq\mu_1<\mu_2\leq k_\kappa}\begin{array}{|c||c|c|c|c|c||c||c|c|c|}\hline\ldots&\multicolumn{5}{c||}{\kappa}&\ldots&\multicolumn{3}{c|}{r+1}\\\hline\ldots&\ldots&&O&&\ldots&\ldots&&\ldots&O\\\hline\ldots&\ldots&\alpha_{\mu_1-1,\kappa}&(i_{\mu_2 \kappa},j_{\mu_1\kappa})&\alpha_{\mu_2+1,\kappa}&\ldots&\ldots&\alpha_{\mu_1+1,\kappa}&\ldots&(i_{\mu_1\kappa},j_{\mu_2,\kappa})\\\hline\ldots&\ldots&&1&&\ldots&\ldots&&\ldots&1\\\hline\end{array}$};

\end{tikzpicture}
\end{figure}
We need again to be careful when we have for instance $\epsilon_{\mu_1\kappa}=0$ and $\epsilon_{\mu_2\kappa}=1$. If we neglect, as always, the problem with the fact of knowing whether or not $s_{\mu_2\kappa}=s_{\mu1+1,\kappa}$, as it still vanishes in the $m\rightarrow\infty$ limit, we obtain :
\begin{eqnarray*}
\ldots&&([O]^{\epsilon_{\mu_1-1\kappa}}_{i_{\mu_1-1\kappa}j_{\mu_1-1\kappa}})_{s_{\mu_1-1\kappa}s_{\mu_1\kappa}}([O]_{i_{\mu_1\kappa}j_{\mu_2\kappa}})_{s_{\mu_1\kappa}s_{\mu_2\kappa}}\ldots([O]^*_{i_{\mu_2\kappa}j_{\mu_1\kappa}})_{s_{\mu_1+1\kappa}s_{\mu_2+1\kappa}}\ldots
\end{eqnarray*}
This can be rewritten:
\begin{eqnarray*}
\ldots&&([O]^{\epsilon_{\mu_1-1\kappa}}_{i_{\mu_1-1\kappa}j_{\mu_1-1\kappa}})_{s_{\mu_1-1\kappa}s_{\mu_1\kappa}}([O]_{i_{\mu_1\kappa}j_{\mu_2\kappa}})_{s_{\mu_1\kappa}s_{\mu_2\kappa}}([O]^{1-\epsilon_{\mu_2-1\kappa}}_{i_{\mu_2-1\kappa}j_{\mu_2-1\kappa}})_{s_{\mu_2\kappa}s_{\mu_2-1\kappa}}\ldots([O]^{1-\epsilon_{\mu_1+1\kappa}}_{i_{\mu_1+1\kappa}j_{\mu_1+1\kappa}})_{s_{\mu_1+2\kappa}s_{\mu_1+1\kappa}}\\
&&\times([O]^*_{i_{\mu_2\kappa}j_{\mu_1\kappa}})_{s_{\mu_1+1\kappa}s_{\mu_2+1\kappa}}\ldots
\end{eqnarray*}
So, this means that we keep only one trace, and because of the factor $1/m$, we get $0$ at the limit, for $\mathcal U^{(2)}_\infty$. For reasons of symmetry, we also obtain zero when for the remaining case.\\
We observe that in the case where $\epsilon_{\mu_1\kappa}=\epsilon_{\mu_2\kappa}$, $U_\infty^{(2)}$ acts like $\hat{\mathcal U}$. This identity does not hold anymore when the $\epsilon$'s differ.

\subsubsection{The operator $\mathcal U^{(1)}$:} Let us first observe that this term is zero whenever $j_{\mu_1\kappa}\neq j_{\mu_2\kappa}$. Moreover, we notice that the term also vanishes in the limit when $\epsilon_{\mu_1\kappa}=\epsilon_{\mu_2\kappa}$, due to an argument similar to what we did with $U^{(2)}_\infty$ when the $\epsilon$ were different. The term becomes non-trivial when the $\epsilon$'s differ, with the creation of a trace. Indeed, computing it, it is easy to see that we have the same thing as what we obtained with $\hat{\mathcal U}$.\\
So, we notice that $\mathcal U^{(1)}+\mathcal U^{(2)}$ correspond in the large $m$ limit exactly to $\hat{\mathcal U}$, with exactly one of both operators being nonzero according  to whether or not the $\epsilon$ differ.

\subsubsection{The operator $\mathcal B^{(1)}$:} Let us first observe again that this term is zero whenever $j_{\mu_1\kappa_1}\neq j_{\mu_2\kappa_2}$. Moreover, we observe that in all the cases of the different values that can be taken by the $\epsilon$'s, the result will be to merge the traces $\kappa_1$ and $\kappa_2$. Therefore, this results is in a factor $1/m^2$, one due to the quadratic variation and another as remnant of the merged trace. So, $\mathcal B_\infty^{(1)}=0$.

\subsubsection{The operator $\mathcal B^{(2)}$:} we need here again to make a difference whether or not $j_{\mu_1\kappa_1}\neq j_{\mu_2\kappa_2}$. When we have $j_{\mu_1\kappa_1}\neq j_{\mu_2\kappa_2}$, there is no problem. When it is not the case, we can reason as in the case of $\mathcal U^{(2)}$ to see that the difference vanishes in the limit $m\rightarrow\infty$. So, we can omitt it. Then, in all four cases (according to the values of the $\epsilon$'s), the traces $\kappa_1$ and $\kappa_2$ merge. Therefore, this term is in $O(1/m^2)$ and, thus, $\mathcal B^{(2)}_\infty=0$.

\subsubsection{Convergence:} We have seen that, at the limit $m\rightarrow\infty$, we obtain the equation $T'=\mathcal D_\infty T+(\mathcal U^{(1)}+\mathcal U^{(2)})T$, which is the same as the equation for the functions $\hat T$ of the quantum stochastic process. Thus, we have the convergence of the marginals of the Brownian motion on $O(nm)$ towards the quantum process $\mathfrak j_t$ in moments. \\
We need now to conclude as to the convergence of the processes, not only of the marginals. This means that we want to prove that $\mathbb{E}(tr_m([O]^{\epsilon_1}_{i_1j_1}(t_1)\ldots [O]^{\epsilon_k}_{i_kj_k}(t_k)))$ converges towards $\phi(\qU^{\epsilon_1}_{i_1j_1}(t_1)\ldots\qU^{\epsilon_k}_{i_kj_k}(t_k))$. The result is already proven when all the $t_l$'s are equal. For the general case, we allow for more general functions. For instance,
$$T=\begin{array}{|c||c|c|c||c|}\hline 1&\multicolumn{3}{c||}{2}&3\\\hline O&m_1&O&m_2&O\\\hline(1,2)&&(2,2)&&(1,1)\\\hline 0&&0&&0\\\hline\end{array}$$
represents $\mathbb{E}\left(tr_m([O]_{1,2}(t))tr_m(m_1[O]_{2,2}(t)m_2)tr_m([O]_{1,1}(t))\right)$, where $m_1$ and $m_2$ are $m\times m$ matrices that are $\sigma(j_m(s),s<t)$-measurable. We then proceed as in \cite[Section 4.4]{Ulrich2015}. These generalized functions $T$ satisfy, in the $m\rightarrow\infty$ limit, the same equations as their generalized quantum counterparts. We then conclude with a recurrence:\\
If $\mathbb{E}(tr_m([O]_{i_1j_1}(t_1)\ldots [O]_{i_kj_k}(t_k)))$ verifies that $t_1=\ldots=t_k$, then it converges towards  $\phi(\qU_{i_1j_1}(t_1)\ldots\qU_{i_kj_k}(t_k))$ by convergence of the marginals.\\
Let us assume that we have shown the convergence of $\mathbb{E}(tr_m([O]_{i_1j_1}(t_1)\ldots [O]_{i_kj_k}(t_k)))$ towards $\phi(\qU_{i_1j_1}(t_1)\ldots\qU_{i_kj_k}(t_k))$ whenever $Card\{t_1,\ldots,t_k\}\leq z$, for a certain integer $z$. Let us now suppose that we have a certain  $\mathbb{E}(tr_m([O]_{i_1j_1}(t_1)\ldots [O]_{i_kj_k}(t_k)))$ with $Card\{t_1,\ldots,t_k\}=z+1$. We can assume, e.g., that the maximum of the $t_l$'s is $t_k$. Then, we can write this expectation as $\mathbb{E}(tr_m(m_1\ldots [O]_{i_kj_k}(t_k)))$, where we replace all the blocks identified in times strictly lesser than $t_k$ by a sign indicating that they are matrices that are $\sigma(j_m(s),s<t_k)$-measurable. By recurrence hypothesis, $([O]_{i,j}(t_l))_{1\leq i,j\leq n,1\leq l\leq k-1 }$ converges in moments towards $(\qU_{i,j}(t_l))_{1\leq i,j\leq n,1\leq l\leq k-1 }$. Now, the function $T$ corresponding to  $\mathbb{E}(tr_m(m_1\ldots [O]_{i_kj_k}(t_k)))$ is entirely characterized by the system of differential equations in which it is involved, as well as by the relationships between the $m_l$'s. The differential equations are still of the same type as before, because the $m_l$'s can be treated as constants. Because by recurrence hypothesis, the matrices $m_l$ do converge towards their quantum counterparts $\tilde m_l$, the relationships existing between them \''converge\'' towards the relationships existing between their quantum counterparts. Thus, $\mathbb{E}(tr_m(m_1\ldots [O]_{i_kj_k}(t_k)))$ converges towards $\phi(\tilde m_1\ldots \qU_{i_kj_k}(t_k))$. This proves the convergence in moments of Theorem \ref{CVortho}.\\

\section{The convergence of the Brownian motion on the symplectic group}
\subsection{The symplectic group and the symplectic dual group}
The symplectic group $Sp(n)$ can be seen in two ways, one of them involving the division ring of quaternions $\mathbb H$. So let us say a few words about the elements of quaternionic linear algebra that we will need. A good and clear reference, on which we rely, is \cite{DV}. We recall that $\HH$ is a division ring, meaning  that it is basically a field but without commutativity. It can be seen as a $\mathbb R$-vector space generated by the family $(1,\imaginaire,\jj,\kk)$ where $1$ is the unit of $\mathbb R$ and $\imaginaire$, $\jj$ and $\kk$ are such that:
$$\imaginaire^2=\jj^2=\kk^2=-1\text{ and }\imaginaire\jj=\kk\text{, }\jj\kk=\imaginaire\text{, }\kk\imaginaire=\jj\text{, }\jj\imaginaire=-\kk\text{, }\kk\jj=-i,\text{, }\imaginaire\kk=-\jj$$
and endowed with a multiplication that generalizes in an obvious way these relations, such that we obtain an algebra.\\
Of course, $\HH$ can be seen either as a $4$-dimensional real vector space, or as a $2$-dimension complex vector space. Especially, as a complex vector space, $\HH$ is generated by $(1,\jj)$. Indeed, a quaternion $q=a+b\imaginaire+c\jj+d\kk$ can be rewritten as $q=(a+b\imaginaire)+\jj(c-d\imaginaire)$. \\There is also an involutive map called conjugation and defined by:
$$\overline{a+b\imaginaire+c\jj+d\kk}=a-b\imaginaire-c\jj-d\kk$$
In particular, for any quaternion $q$, the quantity $q\bar q$ is a positive real number called the modulus of $q$, similar to the modulus of a complex number, and corresponding to the Euclidean distance to the origin from the point of coordinates $(a,b,c,d)$ in the standard real vector space $\mathbb R^4$.\\
We will denote by $\mathcal M_p(\HH)$ the set of $p\times p$ matrices with coefficients in $\HH$. Let $M$ be such a matrix. Because any quaternion $q$ can be rewritten as $q=\alpha+\jj\beta$ with $\alpha$ and $\beta$ complex numbers, it follows that we can write $M$ in the form of $M=Z+\jj W$. Then, we have a map:
$$\begin{array}{cc}\eta:&\mathcal M_p(\HH)\rightarrow\mathcal M_{2p}(\mathbb C)\\&M\mapsto\begin{bmatrix}Z&-\bar W\\W&\bar Z\end{bmatrix}\end{array}$$
This map is a homomorphism of $\mathbb C$-algebras. This basically means that we can adopt one of two equivalent viewpoints: either we see the matrices as being on $\HH$, or we see them as being on $\mathbb C$, but with dimensions doubled. \\
So, from the quaternionic viewpoint, we define the symplectic group $Sp(n)$ as the set of matrices $M\in\mathcal M_n(\HH)$ such that $\,^t\bar M M=I_n$, ie it is the equivalent of the unitary group for the division ring $\HH$. Viewed through $\eta$ one can also say that $Sp(n)$ is the set of matrices $M\in \mathcal M_{2n}(\mathbb C)$ such that $\,^t\bar MM=I_{2n}$ (where the conjugation is the complex one and not the quaternionic one) and $\,^tMJ_nM=J_n$ where $J_n=\begin{bmatrix}0&I_n\\-I_n&0\end{bmatrix}$. This is proven in \cite[(0.4c)]{DV}.\\
How do we define a dual symplectic group $Sp\langle n\rangle$? We must again choose if we take the complex or the quaternionic viewpoint. Voiculescu took the complex one in \cite{Voiculescu87}. Following him, we define $Sp_{\mathbb C}\langle n\rangle$ as the complex algebra $S_n$, endowed with an antilinear involutive automorphism that we will denote by $\bar{\,}$, and generated by $u_{ij},1\leq i,j\leq 2n$ such that the matrix of generators $U$ verifies $U^*U=I_{2n}\,(1)$ and $\bar{U}^*J_nU=J_n\,(2)$, and with the usual comultiplication and coünit.
\begin{rmq}
We need to have an antilinear involution on $S_n$ because $\bar{U}^*$ plays morally the role of $\,^tU$. A complex algebra equipped with such an antilinear involution is what Voiculescu calls a \''real\'' algebra. The idea behind is, I think, that though the algebra is complex, such an antilinear involutive automorphism is an automorphism on a real algebra. Nevertheless, I think that it would be in order to find an other terminology here, especially because the term \''real'\'' is impossible to distinguish from the term real when speaking. I suggest the term realifiable, to emphasize that such an algebra is on the complex, but could very well be said to be on the real numbers.
\end{rmq}
But we could also decide to define a dual symplectic group through the quaternionic viewpoint. Let us define $Sp_{\HH}\langle n\rangle$ as the dual group associated to the $*$-algebra $S^\prime_n$ generated by the $u_{ij},1\leq i,j\leq n$ such that $U^*U=I_n$, but this time the algebra $S^\prime_n$ is taken over $\HH$. A $*$-algebra over the quaternions does not give much more difficulty. We need only to settle for, say, left vector spaces, because $\HH$ is noncommutative.\\

We will in the sequel adopt the quaternionic viewpoint, and thus, $Sp\langle n\rangle$ will always denote $Sp_{\HH}\langle n\rangle$ from now on.\\
What is the stochastic equation of the Brownian motion on $Sp(n)$ seen from the quaternionic viewpoint? The Lie algebra is composed of those matrices $M\in\mathcal M_{nm}(\HH)$ such that $H^*+H=0$. We therefore define a stochastic process $(H(t))_{t\geq0}$ of quaternionic matrices such that:
\begin{enumerate}
\item $M(t)$ is a $nm\times nm$ matrix.
\item $H_{ij}(t)=a_{i,j,t}+\imaginaire b_{ijt}+\jj c_{ijt}+\kk d_{ijt}$ with $a_{ijt}, b_{ijt}, c_{ijt}, d_{ijt}$ independent Gaussian random variables $N(0,\frac{t}{4m}$, for $1\leq i<j\leq nm$.
\item $H_{ij}(t)=-\bar H_{ji}(t)$ for any $i<j$.
\item $H_{ii}(t)=\imaginaire b_{it}+\jj c_{it}+d_{it}\kk$ with $b_{it}, c_{it}, d_{it}$ being independent Gaussian variables $N(0,\frac t{3m})$, for any $1\leq i\leq n$.
\item The family of random variables $H_{ij}$ is independent for $1\leq i\leq j\leq n$.
\item The process $(H(t))_t$ is stationary and with independent increments. 
\end{enumerate}
Let us compute the quadratic variation. We get:
\begin{itemize}
\item[$\bullet$] $d[H_{ab},H_{ab}]=\frac{dt}{4m}-\frac{dt}{4m}-\frac{dt}{4m}-\frac{dt}{4m}=-\frac{dt}{2m}$ for $a\neq b$ (a)
\item[$\bullet$] $d[H_{ab},H_{ba}]=-\frac{dt}{m}$ for $a\neq b$ (b)
\item[$\bullet$] $d[H_{aa},H_{aa}]=-\frac{dt}{m}$ for any $a$ (c)
\item[$\bullet$] $0$ in all other cases
\end{itemize}
Let $a$ and $b$ be two indices. We thus have $(H(t). H(t))_{ab}=\sum_{l=1}^{mn}[H_{al},H_{lb}]=\delta_{ab}.(-\sum_l\frac{dt}{m})+\delta_{ab}.(-\frac{dt}{2m})=-\delta_{ab}(n-\frac1{2m})$
This means that the equation of the Brownian motion $S_t$ on $Sp(nm)$ is:
$$dS_t=S_tdH_t-\frac12(n+\frac1{2m})S_tdt$$
Of course, this yields a tensor-independent Brownian motion on $Sp\langle n\rangle$:
\begin{eqnarray*}
js_m(t):&Sp\langle n\rangle\rightarrow L^{\infty-}(\Omega)\\&u_{ij}\mapsto [S_t]_{ij}
\end{eqnarray*}
When computing the quadratic variations of $S(t)$ we get:
$$d[S_{ab},S_{cd}]=-\underbrace{\delta_{bd}\sum_{s\neq b}S_{as}S_{cs}\frac{dt}{2m}}_{(1)}-\underbrace{\delta_{bd}S_{ab}S_{cb}\frac{dt}{m}}_{(3)}-\underbrace{(1-\delta_{bd})S_{ad}S_{cb}\frac{dt}{m}}_{(2)}$$
We remark that the $(1)$-part is the one due to the quadratic variation (a), the $(2)$ part is due to (c), and the $(3)$ part is due to (c).\\ 
So, similarly to the orthogonal case, we end up with seven operators: $\mathcal D$, $\mathcal U^{(1)}$, $\mathcal U^{(2)}$, $\mathcal U^{(3)}$, $\mathcal B^{(1)}$, $\mathcal B^{(2)}$ and $\mathcal B^{(3)}$. Again, we will be interested only in the $m\rightarrow\infty$ limit of these operators. Let us already examine what happens when there are no $*$. The $(2)$-part is the same as in the orthogonal case. This means that, if $\epsilon_{\mu_1\kappa}=\epsilon_{\mu_2\kappa}=0$ (resp.  $\epsilon_{\mu_1\kappa_1}=\epsilon_{\mu_2\kappa_2}=0$), then  $\mathcal U_ \infty^{(2)}=\hat{\mathcal U}$ (resp.  $\mathcal B_\infty^{(2)}=0$). The $(1)$-part has an additionnal $1/2$ factor, but because these terms vanish, we still have $\mathcal U_\infty^{(1)}=\mathcal B_\infty^{(1)}=0$. Finally, $\mathcal D_\infty$ is clearly the same as $\hat{\mathcal D}$. Indeed, this operator emerges from the $dt$-part of the stochastic equation and in the large $m$-limit, this part reduces to the usual $-\frac n2S_tdt$.\\
What happens if we take $*$ into account? We observe that $\left(\left[S\right]^*_{ij}\right)_{ab}=\left(\left[\bar S\right]_{ij}\right)_{ba}$. Note that this formula was also right in the orthogonal case, but since $\bar O=O$ in this latter case, we ommitted the conjugation.\\
The stochastic equation verifies by $\bar{S}(t)$ is (using the fact that $\bar H_t=-\,^tH_t$):
$$d\bar S_t=-\bar S_td\,^tH_t-\frac12(n+\frac1{2m})\bar S_tdt$$
Thus, we obtain:
$$d[\bar S_{ab},\bar S_{cd}]=-\underbrace{\delta_{bd}\sum_{s\neq b}\bar S_{as}\bar S_{cs}\frac{dt}{2m}}_{(1)}-\underbrace{\delta_{bd}\bar S_{ab}\bar S_{cb}\frac{dt}{m}}_{(3)}-\underbrace{(1-\delta_{bd})\bar S_{ad}\bar S_{cb}\frac{dt}{m}}_{(2)}$$
This means that, if $\epsilon_{\mu_1\kappa}=\epsilon_{\mu_2\kappa}=1$ (resp. $\epsilon_{\mu_1\kappa_1}=\epsilon_{\mu_2\kappa_2}=1$), we obtain the same pattern as for the case where both $\epsilon$'s are equal to zero: $\mathcal U_ \infty^{(2)}=\hat{\mathcal U}$, $\mathcal D_\infty=\hat{\mathcal D}$ and the other operator vanish.\\
We also have:
$$d[\bar S_{ab}, S_{cd}]=-\underbrace{\delta_{bd}\sum_{s\neq b}\bar S_{as} S_{cs}\frac{dt}{2m}}_{(2)}-\underbrace{\delta_{bd}\bar S_{ab} S_{cb}\frac{dt}{m}}_{(3)}-\underbrace{(1-\delta_{bd})\bar S_{ad} S_{cb}\frac{dt}{m}}_{(1)}$$
Let us observe that we have respected the fact that $(1)$ refers to what is obtained from the quadratic variation of $H_t$ of the kind (a), and $(2)$ comes from the quadratic variation of the kind (b). But because we have a $\,^tH_t$ coming into play, the sum over $s$ is now a $(2)$-part.\\
We can then reason as with the orthogonal case and the computations will show that we have $\mathcal U_ \infty^{(1)}=\hat{\mathcal U}$, $\mathcal D_\infty=\hat{\mathcal D}$ and all the other operators vanish.\\
Of course, we can do the same for $d[ S_{ab}, \bar S_{cd}]$ and we still obtain the same thing.\\
So, it is easy to see from that, using the same reasoning with the table-functions and the systems of differential equations, that the marginals of $(js_m(t))_t$ tend towards the marginals of $(\mathfrak j_t)_t$ when $m$ tends towards infinity. We can then apply the same reasoning with the conditional expectation to conclude that $(js_m(t))_t$ converges in $*$-moments towards $(\mathfrak j_t)_t$:
\begin{thm}
When $m$ tends towars infinity, the Brownian motion on the classical symplectic group $Sp(nm)$ tends towards the quantum Lévy process $\mathfrak j$.
\end{thm}

\section*{Conclusion}
We have completed the analysis begun in \cite{Ulrich2015} of the blockwise convergence of the Brownian motion on the classical Lie groups $U(nm)$, $O(nm)$ and $Sp(nm)$. It appears that the Lévy process exhibited in \cite{Ulrich2015,CU2016} plays a special role, because all these Brownian motions converge towards him. In some sense, this seems to mean that the unitary dual group has also a special role, or in other terms, that the dual group taken over the field of complex numbers is privileged over the dual group taken over other division rings.\\
Of course, \cite{Voiculescu87} gives many other examples of dual groups, and they should be studied too in the future. Moreover, it would be interesting to find the relationship, if any, exist between $Sp_{\mathbb C}\langle n\rangle$ and $Sp_\HH\langle n\rangle$. I would conjecture that they ar $\mathbb C$-isomorphic, but was unable to find the isomorphism.\\
As can be seen, much work still needs to be done in the realm of dual groups.
\bibliographystyle{plain}
\bibliography{Convergence_on_orthogonal_and_symplectic_groups}
\nocite{*}

\end{document}